\documentclass[letterpaper, 10 pt, conference]{ieeeconf}
\IEEEoverridecommandlockouts 
\usepackage{amssymb,amsmath,color} 
\usepackage{microtype}
\usepackage{bm}
\usepackage{graphicx} 
\usepackage{subcaption}
\usepackage{epsfig}
\usepackage{algorithm}
\usepackage[noend]{algpseudocode}

\makeatletter
\newcommand{\StatexIndent}[1][3]{%
  \setlength\@tempdima{\algorithmicindent}%
  \Statex\hskip\dimexpr#1\@tempdima\relax}
\makeatother

\algnewcommand\algorithmicinput{\textbf{Initialization:}}
\algnewcommand\init{\item[\algorithmicinput]}

\algnewcommand\algorithmicawake{\textsf{\textit{{AWAKE}}}}
\algnewcommand\awake{\item[\algorithmicawake]}

\algnewcommand\algorithmicidle{\textsf{\textit{{IDLE}}}}
\algnewcommand\idle{\item[\algorithmicidle]}

\newcommand{\broad}{{\small \textbf{BROADCAST }}}

\newcommand{\stopp}{{\small \textbf{STOP }}}

\newcommand{\iidle}{\textsf{\textit{{IDLE }}}}
\newcommand{\aawake}{\textsf{\textit{{AWAKE }}}}

\renewcommand{\natural}{{\mathbb{N}}}

\newcommand{\until}[1]{\{1,\ldots,#1\}}

\newcommand{\EE}{\mathcal{E}}
\newcommand{\GG}{\mathcal{G}} 
 
\newcommand{\NN}{\mathcal{N}}

\newcommand{\VV}{\mathcal{V}}

\newcommand{\st}{\text{subject to }}

\newcommand{\m}{\mathop{\rm minimize}}

\newcommand{\R}{\mathbb{R}}
\newcommand{\p}{\mathbf{P}}
\newcommand{\x}{\mathbf{x}}
\newcommand{\y}{\mathbf{y}}
\newcommand{\LL}{\mathbf{\Lambda}}
\newcommand{\rhoE}{\rho_{E_i}}
\newcommand{\rhoI}{\rho_{I_i}}
\newcommand{\rhoij}{\rho_{ij}}
\newcommand{\tLag}{\tilde{\mathcal{L}}}
\newcommand{\Lag}{\mathcal{L}}

\newtheorem{theorem}{Theorem}[section]
\newtheorem{proposition}[theorem]{Proposition}
\newtheorem{corollary}[theorem]{Corollary}
 \newtheorem{lemma}[theorem]{Lemma}

\newtheorem{assumption}[theorem]{Assumption}
 
\def\QEDclosed{\mbox{\rule[0pt]{1.3ex}{1.3ex}}} 

\def\QED{\QEDclosed} 

\newcommand\oprocendsymbol{\hbox{$\square$}}
\newcommand\oprocend{\relax\ifmmode\else\unskip\hfill\fi\oprocendsymbol}

\renewcommand{\theenumi}{(\roman{enumi})}

\usepackage{tikz}
\usetikzlibrary{shapes,arrows}

\graphicspath{{figs/}}

\newcommand{\AG}[1]{{\color{black} #1}}

\begin{document}

\title{\LARGE \bf
	Asynchronous Distributed Method of Multipliers\\ for Constrained Nonconvex Optimization
}

\author{Francesco Farina$^1$, Andrea Garulli$^1$, Antonio Giannitrapani$^1$, Giuseppe Notarstefano$^2$
	\thanks{$^1$F. Farina, A. Garulli and A. Giannitrapani are with the Dipartimento di Ingegneria dell'Informazione e Scienze Matematiche, Universit{\`a} di Siena, Siena, Italy.}
	\thanks{$^2$G. Notarstefano is with the Department of Engineering,
          Universit{\`a} del Salento, Lecce, Italy.}
\thanks{This result is part of a
    project that has received funding from the European Research Council (ERC)
    under the European Union's Horizon 2020 research and innovation programme
    (grant agreement No 638992 - OPT4SMART).}
}

\maketitle
\thispagestyle{empty}
\pagestyle{empty}


\begin{abstract}
  This paper addresses a class of constrained optimization problems over
  networks in which local cost functions and constraints can be nonconvex. We
  propose an asynchronous distributed optimization algorithm, relying on the
  centralized Method of Multipliers, in which each node wakes up in an
  uncoordinated fashion and performs either a descent step on a local Augmented
  Lagrangian or an ascent step on the local multiplier vector. These two phases
  are regulated by a distributed logic-AND, which allows nodes to understand
  when the descent on the (whole) Augmented Lagrangian is sufficiently small. We
  show that this distributed algorithm is equivalent to a block coordinate
  descent algorithm for the minimization of the Augmented
  Lagrangian followed by an update of the whole multiplier vector. Thus, the proposed
  algorithm inherits the convergence properties of the Method of Multipliers.
\end{abstract}



\section{Introduction}
\label{sec:introduction}
In several cyber-physical network contexts, ranging from control to estimation
and learning, nonconvex optimization problems frequently arise. 
In these contexts, typically each device knows only a portion of the whole
objective function and a subset of the constraints, so that, to avoid the
presence of a central coordinator, distributed algorithms are needed.

Distributed optimization methods relevant for our paper can be divided in two 
groups: those handling nonconvex functions but not local constraints, and
methods handling local constraints, typically designed for convex problems. 
%
Distributed algorithms for nonconvex optimization have started to appear in the
literature only recently. Regarding constrained problems, most of the proposed
methods usually do not deal with local
constraints. In~\cite{bianchi2013convergence} a stochastic gradient method is
proposed,
while in~\cite{wai2016projection} a decentralized
Frank-Wolfe method 
is presented. In~\cite{di2016next,sun2016distributed} the
authors propose distributed algorithms, 
based on the idea of tracking the whole function
gradient and performing successive convex approximations of the nonconvex cost
function. 
%
A perturbed push-sum algorithm for minimizing the sum of nonconvex functions is presented
in~\cite{tatarenko2017non}.

Regarding distributed optimization algorithms handling local constraints,  
in~\cite{lee2013distributed} the authors propose a distributed random
projection algorithm, while a proximal based algorithm is presented
in~\cite{margellos2017distributed}.  In~\cite{necoara2013random} the author
proposes randomized block-coordinate descent methods.
As for
Lagrangian based distributed algorithms, in~\cite{wei20131} an asynchronous
distributed ADMM is proposed for a separable, constrained optimization
problem. Other ADMM based approaches are presented
in~\cite{iutzeler2016explicit,bianchi2014stochastic,bianchi2016coordinate},
while an asynchronous proximal dual algorithm is proposed in
\cite{notarnicola2017asynchronous}.
Finally, a distributed algorithm for a structured class of nonconvex
optimization problems with local constrains is presented
in~\cite{notarnicola2016randomized}. 

The main contribution of this paper is a \emph{fully} distributed optimization
algorithm called ASYnchronous Method of Multipliers (ASYMM), addressing
constrained optimization problems over networks 
in which both local cost
functions and local constraints may be nonconvex. 
%
%
ASYMM is designed for asynchronous networks. It features two types of local updates at each node, a primal and a multiplier one, which
are regulated by an asynchronous distributed logic-AND algorithm. 
When awake, nodes start performing a descent step on a local Augmented Lagrangian.
By means of a distributed logic-AND algorithm, they realize when all
of them have reached a given tolerance on their local gradient. Thus, still asynchronously, they start
performing the multiplier updates. 
An interesting feature of ASYMM is that a
node 
just needs to receive all neighboring multipliers to start again its primal
descent.
The analysis shows that ASYMM implements a
suitable inexact version of the Method of Multipliers 
and hence it
inherits all the convergence properties of the centralized method (see~\cite{rockafellar1974augmented, bertsekas2014constrained}). 
%
It is worth mentioning
that our algorithm implicitly leverages on some results given in~\cite{matei2014extension,matei2017nonlinear} to handle non-regularity of local minima (due to the presence of local copies of the common decision variable). In this respect, ASYMM presents three main novelties with respect to the
distributed algorithm proposed in
\cite{matei2014extension,matei2017nonlinear}: the network model is asynchronous;
primal minimization is performed approximately; switching from a primal to a
multiplier update is performed in a distributed way. 

Due to space limitations all the proofs of the presented results are omitted and will be provided in a forthcoming document.

{\bf Organization}: 
In Section~\ref{sec:setup} we present the distributed optimization set-up. The
proposed distributed algorithm is presented in Section~\ref{sec:algorithm} and
analyzed in Section~\ref{sec:analysis}. Finally, a numerical application is
presented in Section~\ref{sec:numerical}.

{\bf Notation and definitions}: Given a matrix $A\in\R^{n\times m}$ we denote by
$A[i,j]$ the $(i,j)$-th element of $A$, $A[:,i]$ its $i$-th column and $A[i,:]$
its $i$-th row. We write $A\left[i,:\right]=b$ to assign the value $b$ to all
the elements in the $i$-th row of $A$. Given two vectors $A,B\in\R^{n\times 1}$, we write $A>c$ if all elements of $A$ 
are greater than $c$ and $A>B$ if $A[i]>B[i]$ for all $i$. We denote by
$\mathbf{0}_{n\times m}$ the $n\times m$ zero matrix. 
If $J=\left\lbrace j_1,...,j_m \right\rbrace$ is a set of indexes, we denote by
$[z_j]_{j\in J}$ the vector $[z_{j_1},...,z_{j_m}]$.\\
%
A function $\Psi(x)$
has \emph{Lipschitz continuous gradient} if there exists a constant $L$ such
that $\|\nabla\Psi(x) - \nabla\Psi(y)\|\leq L\|x-y\|$ for all $x,y$.  It is
\emph{$\sigma$-strongly convex} if
$(\nabla\Psi(x) - \nabla\Psi(y))^\top(x-y)\geq \sigma\|x-y\|^2$. 
Let $\x=[x_1^\top,...,x_N^\top]^\top$, with $x_i\in\R^{n_i}$ and
$\sum_{i=1}^N n_i = m$, and let $U^{m\times m}$ be a column \AG{partition} of the
identity matrix such that $\x = \sum_{i=1}^N U_i x_i$ and $x_i=U_i^\top\x$.
The function $\Phi(\x)$ has \emph{block component-wise Lipschitz continuous gradient} if
there are constants $L_i\geq 0$ such that
$ \|\nabla_{x_i}\Phi(\x+U_i s_i) - \nabla_{x_i}\Phi(\x)\|\leq L_i\|s_i\| $
for all $\x\in \R^{n}$ and $s_i\in\R^{n_i}$.

\section{Set-up and Preliminaries}

\subsection{Distributed Optimization Set-up}
\label{sec:setup}
Consider the following optimization problem
\begin{equation}\label{pb:problem}
\begin{aligned}
& \m_x
& & \sum_{i=1}^N f_i(x)\\
& \st
& & h_i(x)=0, & i=1,...,N,\\
& & & g_i(x)\leq 0, & i=1,...,N,
\end{aligned}
\end{equation}
with $f_i, h_i, g_i:\R^n\to\R$ satisfying the following assumptions.

\begin{assumption}\label{asm:C2_functions}
  Functions $f_i,h_i,g_i\in C^2$, $i\in\{1,...,N\}$ have bounded Hessian and
  Lipschitz continuous gradients. Moreover, problem~\eqref{pb:problem} has at
  least one feasible solution.\oprocend
\end{assumption}

\begin{assumption}\label{asm:sosc}
  Every local minimum of~\eqref{pb:problem} is a regular point and
  satisfies the second order sufficiency conditions. \oprocend
\end{assumption}

Problem~\eqref{pb:problem} is to be solved in a distributed way by a network of $N$ peer
processors without a central coordinator. Each processor has a local memory, a
local computation capability and can exchange information with neighboring
nodes. Moreover, functions $f_i$, $h_i$ and $g_i$ are private to node $i$.
The network is described by a fixed, undirected and connected graph
$\GG=(\VV,\EE)$, where $\VV=\{1,...,N\}$ is the set of nodes and
$\EE\subseteq\{1,...,N\}\times\{1,...,N\}$ is the set of edges.  We denote by
$\NN_i=\{j\in \VV\mid (i,j)\in\EE\}\cup \{i\}$ the set of the neighbors of node
$i$ (including $i$ itself) and by $d_i=|\NN_i|$ the cardinality of
$\NN_i$. Also, we denote by $d_G$ the diameter of $\GG$.

We consider a generalized version of the
asynchronous communication protocol presented in~\cite{notarnicola2017asynchronous}. Each node
has its own concept of time defined by a local timer, which triggers when the node has
awake, independently of the other nodes.  Between two triggering events
each node is in \iidle mode, i.e., it listens for messages from neighboring
nodes and, if needed, updates some local variables.
When a trigger occurs, it switches into \aawake mode, performs
local computations and broadcasts the updated information to neighbors.
Formally, the triggering process is modeled by means of a local clock
$\tau_i\in\R_{>0}$ and a certain waiting time $T_i$. As long as $\tau_i<T_i$ the
node is in \iidle. When $\tau_i =T_i$ the node switches to the \aawake mode
and, after running the local computations, resets $\tau_i = 0$ and selects a new
waiting time $T_i$, according to the following assumption.
\begin{assumption}[Local timers]\label{asm:timing}
  For each node $i$, there exists $\bar{T}_i$ such that $T_i\leq \bar{T}_i$ for
  each awakening cycle.\oprocend
\end{assumption}
\begin{assumption}[No simultaneous awakening]\label{asm:no_awake}
  Only one node can be awake at each time instant. \oprocend
\end{assumption}



\subsection{Equivalent Formulation and Method of Multipliers}
Due to the connectedness of $\GG$, problem~\eqref{pb:problem}, can be rewritten
in the equivalent form
\begin{equation}\label{pb:distributed_problem}
\begin{aligned}
&\m_{x_1,...,x_N}
& &  \sum_{i=1}^N f_i(x_i)\\
& \st
& & x_i=x_j,& \forall(i,j)\in \EE,\\
& & & h_i(x_i)=0, & i\in\VV,\\
& & & g_i(x_i)\leq 0, & i\in\VV.
\end{aligned}
\end{equation}

Next we define the Augmented Lagrangian associated to
problem~\eqref{pb:distributed_problem}.
Let $\nu_{ij}\in\R^n$ and $\rhoij$ be the multiplier and penalty parameter
associated to the equality constraint $x_i = x_j$. Similarly, let $\lambda_i$
and $\rhoI$, respectively $\mu_i$ and $\rhoE$, be the multiplier and penalty
parameter associated to the $i$-th inequality, respectively equality,
constraint.
Moreover, let $\x=[x_1^\top,...,x_N^\top]^\top$, denote by $\p$ the vector
stacking all the penalty parameters, $\nu$, $\lambda$ and $\mu$ the vectors
stacking the corresponding multipliers, and, consistently let
$\LL=[\nu^\top,\lambda^\top,\mu^\top]^\top$.  Thus, the Augmented Lagrangian
associated to~\eqref{pb:distributed_problem} is
\begin{align}
\Lag_{\p}(\x,\LL)=&\sum_{i=1}^N \bigg\lbrace f_i(x_i) \nonumber\\
&+\sum_{j\in \NN_i\backslash i}\left[\nu_{ij}^\top(x_i-x_j)+\frac{\rhoij}{2}\|x_i-x_j\|^2\right]+\nonumber\\
&+\lambda_i h_i(x_i)+\frac{\rhoE}{2}\|h_i(x_i)\|^2+\nonumber\\
&+\frac{1}{2\rhoI}\left(\max\{0,\, \mu_i + \rhoI g_i(x_i) \}^2-\mu_i^2\right)\bigg\rbrace.\label{eq:L}
\end{align}

A powerful method to solve problem~\eqref{pb:distributed_problem} is the well
known Method of Multipliers, which consists of the following steps (see
e.g. \cite{rockafellar1974augmented, bertsekas2014constrained}),
\begin{align}
	\x^{k+1} &= \arg \min_{\x} \Lag_{\p^k}(\x,\LL^k)\label{eq:x_minimization}\\
	\nu_{ij}^{k+1} &= \nu_{ij}^k+\rho_{ij}^k(x_i^{k+1}-x_j^{k+1}), &\forall(i,i)\in\EE,\label{eq:nu}\\
	\lambda_i^{k+1} &= \lambda_i^k + \rhoE^k h_i(x_i^{k+1}),&i\in\VV,\\
	\mu_i^{k+1} &= \max\{0,\, \mu_i^k+\rhoI^k g_i(x_i^{k+1})\},& i\in\VV, \label{eq:mu}
\end{align}
where $\p^{k+1}\geq\p^k\geq...\geq\p^0> 0$. 
Sufficient conditions guaranteeing the convergence of this method to a local
minimum of problem~\eqref{pb:distributed_problem} have been given, e.g., in
\cite{bertsekas2014constrained}.  One of these conditions involves the
regularity of the local minima of the optimization problem. In general, such a
condition is not verified in problem~\eqref{pb:distributed_problem} due to the
constraints $x_i=x_j$. 
In~\cite{matei2014extension,matei2017nonlinear} the results
in~\cite{bertsekas2014constrained} have been extended to deal with the non
regularity of the local minima of problem~\eqref{pb:distributed_problem}.



We want to stress that the Augmented Lagrangian defined in~\eqref{eq:L} is not
separable in the local decision variables $x_i$. Thus, the minimization step
in~\eqref{eq:x_minimization} cannot be performed by independently minimizing
with respect to each variable.


\section{Asynchronous Method of Multipliers Distributed Algorithm}
\label{sec:algorithm}
In this section, the Asynchronous Method of Multipliers (ASYMM) for solving problem~\eqref{pb:distributed_problem} in an asynchronous
and distributed way is presented. We start by introducing a distributed
algorithm allowing nodes in an asynchronous network to agree that all of
them have set a local flag to $1$.

\subsection{Asynchronous distributed logic-AND} 
Each node in the network is assigned a flag
$C_i$ that is initially set to $0$ and is then changed to $1$ in finite time.
We propose an asynchronous distributed logic-AND algorithm, based on the
synchronous algorithm proposed in~\cite{ayken2015diffusion}, for checking if all
the nodes have $C_i=1$.

%
Each node $i$ stores a matrix $S_i\in\{0,1\}^{d_G\times d_i}$ which contains
information about the status of the node itself and its neighbors.  Let
$S_i[l,j|_i]$ denote the element in the $l$-th row and $j|_i$-th column of
$S_i$, where $j|_i$ is the index associated to node $j$ by node $i$.  Then, the
elements $S_i[1,j|_i]$ for $j\in \NN_i\backslash i$ represent the values of the
flags of nodes $j\in\NN_i\backslash i$ and $S_i[1,d_i]$ the one of node $i$
itself. This means that $S_i[1,d_i]=C_i$.
Moreover, for $l=2,\dots,d_G$, the element $S_i[l,j|_i]$,
$j\in \mathcal{N}_i \setminus i$, represents the status of the $(l-1)$-th row of
$S_j$ defined as the product of all its entries. Similarly, $S_i[l,d_i]$
represents the status of the $(l-1)$-th row of $S_i$ and can be computed as
\begin{equation}
S_i[l,d_i]=\prod_{b=1}^{d_i}S_i[l-1,b].
\end{equation}
Hence, $S_i[l,d_i]=1$ if and only if $S_i[l-1,j|_i]=1$ $\forall j\in \NN_i$. A pseudo code of the distributed logic-AND algorithm is reported in Algorithm~\ref{alg:async_AND}. 
Notice, in particular, that node $i$ has to broadcast to all
its neighbors only the last column of $S_i$, i.e. $S_i[:,d_i]$. Moreover it stores only the $d_j$-th column of the matrices $S_j$
of its neighbors $j\in \NN_i\backslash i$.
\begin{algorithm}
	\begin{algorithmic}
          \init $C_i\gets 0$, 
          $S_i \gets \mathbf{0}_{d_G\times d_i}$
		\item[]
		\awake
		\If{$\prod_{l=1}^{d_i}S_i[d_G,l]\neq1$}
		\State $S_i[1,d_i]\gets C_i$
		\State $S_i[l,d_i]\gets\prod_{b=1}^{d_i}S_i[l-1,b]$ for $l=2,...,d_G$
		\State \broad $S_i[:,d_i]$ to all $j\in \NN_i\backslash i$
		\EndIf
		\item[]
		\If{$\prod_{l=1}^{d_i}S_i[d_G,l]=1$}
		\State \stopp and send \stopp signal to all $j \in \NN_i\backslash i$
		\EndIf
              \item[] \idle \If{$S_j[:,d_j]$ received from
                  $j\in \NN_i\backslash i$ and not received a \stopp
                  signal} 
		\State $S_i[l,j|_i]\gets S_j[l,d_j]$ for $l=1,...,d_G$
		\EndIf
		\State \textbf{if} \stopp received, set $S_i[d_G,:]\gets 1$

	\end{algorithmic}
	\caption{Asynchronous distributed logic-AND}\label{alg:async_AND}
\end{algorithm}

It can be seen that a node will stop only when the last row of its stopping
matrix is composed by all $1$s, i.e. when
\begin{equation}\label{eq:stop_S_equation}
\prod_{b=1}^{d_i}S_i[d_G,b]=1.
\end{equation}

The following result states that~\eqref{eq:stop_S_equation} can be satisfied
at some node if and only if $C_i=1$ for all $i$.

\begin{proposition}
\label{prop:logic-AND}
  Let Assumption~\ref{asm:timing} holds. 
  Suppose that when a flag $C_i$ switches from $0$ to $1$, it remains equal to $1$ indefinitely.
  Then, a node $i$ satisfies $\prod_{b=1}^{d_i}S_i[d_G,b]=1$ in finite time if
  and only if at a certain time instant one has $C_j=1$ for all $j\in\VV$. \oprocend
\end{proposition}

\subsection{ASYMM Algorithm}
We now present the proposed ASYMM algorithm for solving
problem~\eqref{pb:distributed_problem} in an asynchronous and distributed way.

Before describing the algorithm, we need to introduce a ``local Augmented
Lagrangian'', whose minimization with respect to the decision variable $x_i$ is
equivalent to the minimization of the entire Augmented Lagrangian
\eqref{eq:L}. 
To this end, define $x_{\NN_i} = [x_{j}]_{j\in \NN_i}$,
$\LL_{\NN_i}=[\lambda_i,\mu_i,[\nu_{ij}]_{j\in \NN_i\backslash i},
[\nu_{ji}]_{j\in \NN_i\backslash i}]$,
and
$\p_{\NN_i}=[\rhoI,\rhoE,[\rhoij]_{j\in \NN_i\backslash i},[\rho_{ji}]_{j\in
  \NN_i\backslash i}]$.
Then we can introduce the \emph{local Augmented Lagrangian}
\begin{align}
\tLag_{\p_{\NN_i}}&(x_{\NN_i},\LL_{\NN_i})=\nonumber\\
&= f_i(x_i)+ \nonumber\\
&+\sum_{j\in \NN_i\backslash i}\left[x_i^\top(\nu_{ij}-\nu_{ji})+\frac{\rhoij+\rho_{ji}}{2}\|x_i-x_j\|^2\right]+\nonumber\\
&+\lambda_i h_i(x_i)+\frac{\rhoE}{2}\|h_i(x_i)\|^2+\nonumber\\
&+\frac{1}{2\rhoI}\left(\max\{0,\, \mu_i + \rhoI g_i(x_i) \}^2-\mu_i^2\right).\label{eq:L_local}
\end{align}

It can be easily verified that, for fixed values of $x_j$, $j\neq i$,
minimizing~\eqref{eq:L} with respect to $x_i$ is equivalent to
minimizing~\eqref{eq:L_local} with respect to $x_i$, i.e.
\begin{equation*}
	\arg \min_{x_i}\Lag_{\p}(\x,\LL)= \arg\min_{x_i}\tLag_{\p_{\NN_i}}(x_{\NN_i},\LL_{\NN_i}).
\end{equation*}
Moreover, the gradients with respect to $x_i$ are equal, i.e.
\begin{equation}\label{eq:gradient_equivalence}
	\nabla_{x_i}\Lag_{\p}(\x,\LL) = \nabla_{x_i}\tLag_{\p_{\NN_i}}(x_{\NN_i},\LL_{\NN_i}).
\end{equation}

We are now ready to introduce ASYMM (which we report in Algorithm~\ref{alg:async}), whose rationale
is the following. When a node wakes up, it performs a gradient descent step on its
local Augmented Lagrangian until every node has reached a suitable accuracy.
This check is performed by nodes themselves in a distributed way. When a node
gets aware of this condition, it performs (only) one ascent step on its local
multiplier vector. Then, it gets back to the primal update when it has received
the updated multipliers from all its neighbors. 

More formally, when node $i$ wakes up, it checks through a flag, called $M_{done}$,
if its multiplier vector and the neighboring ones are up to date. If this is the case (which corresponds to $M_{done}=0$),
it performs one of the following two tasks:
\begin{enumerate}
\item[T1.] \label{task:1} If $\prod_{l=1}^{d_i}S_i[d_G,l]\neq 1$, it performs a
  gradient descent step on its local Augmented Lagrangian and checks if the
  local tolerance on the gradient has been reached (if the latter is true this corresponds to setting
  $C_i\gets 1$ in the distributed logic-AND). Then
  it updates matrix $S_i$, and broadcasts the updated $x_i$ and $S_i[:,d_i]$ to its
  neighbors.
\item[T2.] \label{task:2} If $\prod_{l=1}^{d_i}S_i[d_G,l] = 1$, it performs an
  ascent step on the local multiplier vector and updates the local penalty
  parameters, then sets $M_{done}=1$ and broadcasts the updated multipliers and
  penalty parameters $\nu_{ij}$ and $\rho_{ij}$ (associated to constraints
  $x_i=x_j$) to its neighbors.
\end{enumerate}

When in \iidle, node $i$ continuously listens for messages from its neighbors,
but does not broadcast any information. Received messages may contain either the
local optimization variable and variables for the logic-AND, or multiplier
vectors and penalty parameters. If necessary, node $i$ suitably updates local
variables of the logic-AND or the flag $M_{done}$.  We note that, for node $i$,
sending a new multiplier $\nu_{ij}$ or receiving a new $\nu_{ji}$ corresponds to
sending or receiving a \stopp signal in the asynchronous logic-AND.  Finally,
regarding the rule used for updating the penalty parameters, the heuristics
presented, e.g., in~\cite[Chapters~2 and~3]{bertsekas2014constrained} can be used.


\begin{algorithm}
	\caption{ASYMM}\label{alg:async}
	\begin{algorithmic}
		\init Initialize $x_i$, $\LL_i$, $\bm{\nu}_i$, $\p_i$, $\bm{\rho}_i$, $S_i=\mathbf{0}_{d_G\times d_i}$, $M_{done}$ = 0.
		\item[]
		
		\awake
		\If{$\prod_{l=1}^{d_i}S_i[d_G,l]\neq 1$ \textbf{and} \textbf{not} $M_{done}$}\\
		\vspace{-1ex}
		\State $x_i\gets x_i-\frac{1}{L_i}\nabla_{x_i}\tLag_{\p_{\NN_i}}(x_{\NN_i},\LL_{\NN_i})$\\
		\vspace{-1ex}
		\If{$\|\nabla_{x_i}\tLag_{\p_{\NN_i}}(x_{\NN_i},\LL_{\NN_i})\|\leq\epsilon_i$} $S_i[1,d_i]\gets 1$\\
		\vspace{-1.5ex}
		\EndIf
		\State $S_i[l,d_i]\gets\prod_{b=1}^{d_i}S_i[l-1,b]$ for $l=2,...,d_G$\\
		\vspace{-1.5ex}
		\State \broad $x_i$, $S_i[:,d_i]$ to all $j\in \NN_i\backslash i$
		\EndIf
		\item[]
		\If{$\prod_{l=1}^{d_i}S_i[d_G,l]=1$ \textbf{and} \textbf{not} $M_{done}$}\\
		\vspace{-1ex}
		\State $\nu_{ij}\gets\nu_{ij}+\rho_{ij}(x_i-x_j)$ for $j\in \NN_i\backslash i$\\
		\vspace{-2ex}
		\State $\lambda_{i}\gets \lambda_{i}+\rho_{E_i} h_i(x_i)$\\
		\vspace{-2ex}
		\State $\mu_{i}\gets \max\{0,\, \mu_{i}+\rho_{I_i} g_i(x_i)\}$\\
		\vspace{-1ex}
		\State update $\rhoE$, $\rhoI$, $\rhoij$ for all $j\in \NN_i\backslash i$
		\State $M_{done}$ $\gets$ 1
		\State \broad $\nu_{ij}$, $\rhoij$ to $j \in \NN_i\backslash i$ 
		\EndIf
		\item[]
		
		\idle
		\If{$S_j[:,d_j]$ received from $j\in \NN_i\backslash i$ and not already received some new $\nu_{ji}$} 
		 $S_i[l,j|_i]\gets S_j[l,d_j]$ for $l=1,...,d_G$
		\EndIf
		\State \textbf{if} $\nu_{ji}$ and $\rho_{ji}$ received from $j\in \NN_i\backslash i$ set $S_i\left[d_G,:\right]\gets 1$
		\State \textbf{if} $x_j^{new}$ received from $j\in \NN_i\backslash i$, update $x_j\gets x_j^{new}$
		\If{$M_{done}$ \textbf{and} $\nu_{ji}$ received from all $j\in \NN_i\backslash i$}
		\State $M_{done}$ $\gets$ 0, $S_i\gets\mathbf{0}_{d_G\times d_i}$, update $\epsilon_i$
		\EndIf
	\end{algorithmic}
\end{algorithm}
\section{ASYMM Convergence Analysis}
\label{sec:analysis}
In order to analyze ASYMM, we start by noting that under
Assumption~\ref{asm:timing}, from a global perspective, the local asynchronous
updates result into an algorithmic evolution in which, at each iteration, only
one node wakes up in an essentially cyclic fashion\footnote{Indexes in
  $\until{N}$ are drawn according to an \emph{essentially cyclic} rule if there
  exists $M\geq N$ such that every $i\in\until{N}$ is drawn at least once every
  $M$ extractions.}. Hence, we can associate to each triggering an iteration of
the distributed algorithm. We denote by $t\in\natural$ a discrete,
\emph{universal} time indicating the $t$-th iteration of the algorithm and
define as $i_t\in\VV$ the index of the node triggered at iteration
$t$. 

In the following we show that: (i) there is an \emph{equivalence} relationship
between ASYMM and an inexact Method of Multipliers and (ii) the
\emph{convergence} of ASYMM to a local minimum of problem~\eqref{pb:problem} is
guaranteed under suitable conditions inherited from the (centralized)
optimization literature.

\subsection{Equivalence with an inexact Method of Multipliers}
We consider an inexact Method of Multipliers which consists of solving the $k$-th
instance of the Augmented Lagrangian minimization by means of a block-coordinate
gradient descent algorithm (see, e.g., \cite{wright2015coordinate} for a
survey), which runs for a certain number of iterations $h^k$. A pseudo code of
this inexact Method of Multipliers is given in the following table, where
$i_{h}$ is the index of the block chosen at iteration $h$ and the penalty parameters are updated as in~\cite{bertsekas2014constrained}.

\begin{algorithm}
	\begin{algorithmic}
		\For{$k=0,1,...$}
		\State $\hat{\x}^0=\x^k$
		\For{$h =1,...,h^k$}
		\State $\hat{\x}^{h+1}=\hat{\x}^h-\frac{1}{L_{i_h}} U_{i_h}\nabla_{x_{i_h}}\Lag_{\p}(\hat{\x}^h,\LL^k)$
		\EndFor
		\State $\x^{k+1}=\hat{\x}^{h^k+1}$
		\State $\nu_{ij}^{k+1} = \nu_{ij}^k+\rho_{ij}^k(x_i^{k+1}-x_j^{k+1}),\, \forall(i,i)\in\EE$
		\State $\lambda_i^{k+1} = \lambda_i^k + \rhoE^k h_i(x_i^{k+1}),\,i\in\VV$
		\State $\mu_i^{k+1} = \max\{0,\, \mu_i^k+\rhoI^k g_i(x_i^{k+1})\},\, i\in\VV$
		\EndFor	
	\end{algorithmic}
	\caption{Inexact MM}\label{alg:inexactMM}
\end{algorithm}

We would like to clarify that the ordered sequence of indexes $h$ and $k$ used
in Algorithm~\ref{alg:inexactMM} does not coincide with the sequence in the
universal time $t$ of ASYMM algorithm. We will rather show that a (possibly reordered)
subsequence of iterations (in the universal time $t$) of ASYMM will give rise
to the $h$ and $k$ sequences in Algorithm~\ref{alg:inexactMM}.

Let $t_1,t_2,...$ be a subsequence of $\{t\}$ such that at each $t_\ell$ a
multiplier update (task T2) has been performed by node $i_{t_\ell}$ and let $t_1$
be the time instant of the first multiplier update. Then, the following result
holds.  
\begin{lemma}\label{lemma:l_cyc} %
  Each sequence $(i_{t_{kN+1}},...,i_{t_{(k+1)N}})$, for $k=0,1,...$, is a
  permutation of $\{1,...,N\}$. Moreover, if $\epsilon_i>0$ $\forall i\in\VV$,
  multiplier updates occur infinitely many times. \oprocend
\end{lemma}

Define $\LL_i=[\lambda_i,\mu_i,[\nu_{ij}]_{j\in\NN_i\backslash i}]$ and let
$\tilde{x}_{i}^t$ and $\tilde{\Lambda}_{i}^t$ be the value of the state vector and of the
multiplier vector at node $i$, at iteration $t$, computed according to ASYMM.  Then, the following
Corollary follows immediately from Lemma~\ref{lemma:l_cyc}.

\begin{corollary}
\label{cor:const_lambda}
  For all $\tau\in\{t_{kN+1},...,t_{(k+1)N}\}$, $k=0,1,...$, it holds
  $\tilde{\LL}_{i_\tau}^{t}=\tilde{\LL}_{i_\tau}^{\tau}$ $\forall t\!\in\!\!\{\tau,\tau\!+\!1,\dots,t_{(k+1)N}\}$. \oprocend
\end{corollary}

For all $\tau\in\{t_{kN+1},...,t_{(k+1)N}\}$, $k=0,1,...$, define
\begin{equation*}
	x_{i_\tau}^{k+1}=\tilde{x}_{i_\tau}^\tau,\quad
	\LL_{i_\tau}^{k+1}=\tilde{\LL}_{i_\tau}^\tau.
\end{equation*}
By using Lemma~\ref{lemma:l_cyc} and reordering the indexes $i_{\tau}$, one can define
\begin{align*}
	\x^{k+1}&=\left[\left(x_1^{k+1}\right)^\top,...,\left(x_N^{k+1}\right)^\top\right]^\top,\\
	\LL^{k+1}&=\left[\left(\LL_1^{k+1}\right)^\top,...,\left(\LL_N^{k+1}\right)^\top\right]^\top.
\end{align*}

The next two lemmas show that a local primal (resp. multiplier) update is performed
according to a common multiplier (resp. primal) variable.
\begin{lemma}\label{lemma:l_consistency}
  For all $\tau\in\{t_{kN+1},...,t_{(k+1)N}\}$, $k=0,\!1,...$, every multiplier
  update $\LL_{i_\tau}^{k+1}$ is performed using $\x^{k+1}$. \oprocend
\end{lemma}

\begin{lemma}\label{lemma:x_consistncy}
  Let $\tau_i^k$ be the time instant in which $x_i^{k+1}$ is computed, i.e.,
  $\tau_i^k \in \{t_{kN+1},...,t_{(k+1)N}\}$ such that node $i$ is awake at time
  $\tau_i^k$. Then, for all $t\in\{\tau_i^k,\tau_i^k+1,\dots,\tau_i^{k+1}\}$, every
  descent step on the Augmented Lagrangian with respect to 
  $x_i$ is performed using the multiplier vector $\LL^{k+1}$. \oprocend
\end{lemma}

Next lemma states that every node performs at least one primal update between
the beginning of two consecutive cycles of multiplier updates.
\begin{lemma}\label{lemma:awake_once}
  Between $t_{kN+1}$ and $t_{(k+1)N+1}$ every node performs task $T1$ at least
  once. \oprocend
\end{lemma}

The equivalence of ASYMM and Algorithm~\ref{alg:inexactMM} is stated in the next theorem, whose proof relies on the previous Lemmas.
\begin{theorem}\label{thm:equivalence}
  Let Assumptions~\ref{asm:C2_functions}, \ref{asm:sosc}, \ref{asm:timing} and
  \ref{asm:no_awake} hold.
  Then, ASYMM is equivalent to an instance of Algorithm~\ref{alg:inexactMM} in
  which the selection of nodes $i_h$ satisfies an essentially cyclic rule.
%
Moreover, if in Algorithm~\ref{alg:async}, $\epsilon_i>0$ $\forall i\in\VV$,
  the total number of primal descent steps $h^k$ is finite. 
  \oprocend
\end{theorem}
%
%
 
\subsection{Local Convergence}
Let us first introduce a result that allows us to bound the norm of the gradient of a strongly convex
function with block component-wise Lipschitz continuous gradient, during the
evolution of a block coordinate descent algorithm. Specifically, we relate this bound to given local
tolerances $\epsilon_i$ on the norm of the gradient with respect to block $i$.
\begin{lemma}\label{thm:local_gradient_bound}
  Let $\Phi(y_1, \ldots, y_N)$ be a $\sigma$-strongly convex function with block
  component-wise Lipschitz continuous gradients (with $L_i$ being the Lipschitz
  constant with respect to block $y_i$) in a subset $Y\subseteq\R^n$.
  Let $\{\y^h\}$ be a sequence generated starting according to $\y^{h+1}=\y^h-\frac{1}{L_{i_h}} U_{i_h}\nabla_{y_{i_h}}\Phi(\y^h)$,
  where $\y^0\in Y$ and indexes $i_h\in\until{N}$ are drawn in an essentially cyclic way.
  If, for some $\bar{h}>0$, $\|\nabla_{y_i}\Phi(\y^{\bar{h}})\|\leq \epsilon_i, \; \forall i \in
\until{N}$,
  then
  \begin{equation*}
    \|\nabla_{\y}\Phi(\y^h)\|\leq
\sqrt{\sum_{i=1}^N\left(\frac{L_i\epsilon_i}{\sigma}\right)^2}
  \end{equation*}
  for all $h\geq\bar{h}$. \oprocend
\end{lemma}

In order to show the local convergence, we need an additional assumption.
\begin{assumption}
\label{asm:xk_local}
There exists some $\bar{k}>0$, such that for all $k\geq\bar{k}$, the sequence $\x^k$ generated
by ASYMM belongs to a $\sigma^k$-strongly convex neighborhood of a local minimum
of $\Lag_{\p^k}(\x,\LL^{k})$. \oprocend
\end{assumption}

Assumption~\ref{asm:xk_local} is indeed strong, 
but it is somehow standard in the optimization literature. As pointed out, e.g., in
\cite{bertsekas2014constrained}, while it is not possible to guarantee that such
an assumption holds a priori, in practice it turns out to be 
satisfied after a sufficient number of iterations of the primal minimization and
multiplier/penalties update.
The next Theorem, whose proof relies on Lemma~\ref{thm:local_gradient_bound}, shows that ASYMM guarantees
$\|\nabla_{\x}\Lag_{\p^k}(\x^{k+1},\LL^k)\| \leq \varepsilon^k$, with $\varepsilon^k$
depending on the local thresholds $\epsilon_i$.
\begin{theorem}\label{thm:conv}
 Let Assumptions~\ref{asm:C2_functions}, \ref{asm:sosc}, \ref{asm:timing},
  \ref{asm:no_awake} and~\ref{asm:xk_local} hold. Then, there exists $\bar{k}>0$, such that
  for all $k\geq\bar{k}$, it holds
  $\|\nabla_{\x}\Lag_{\p^k}(\x^{k+1},\LL^k)\| \leq \varepsilon^k$ with
	\begin{equation*}
	\varepsilon^k
	=\sqrt{\sum_{i=1}^N\left(\frac{L_i^k\epsilon_i^k}{\sigma^k}\right)^2},
	\end{equation*} 
  where $\epsilon_i^k$ is the local tolerance set by node $i$ for the primal
  descent related to multiplier $\LL^k$ and $L_i^k$ is the Lipschitz
  constant of $\nabla_{\x}\Lag_{\p^k}(\x,\LL^k)$ with respect to $x_i$. \oprocend
\end{theorem}

We want to remark that the only global parameter $\sigma^k$ appearing in our
analysis does not need to be known by the nodes, because it is not required for the execution of ASYMM. 

We stress that ASYMM is
equivalent to Algorithm~\ref{alg:inexactMM} even without
Assumption~\ref{asm:xk_local}, which is needed to guarantee the local convergence to a
(strict) local minimum of problem~\eqref{pb:problem}. In fact, for $k<\bar{k}$ in Assumption~\ref{asm:xk_local} the Augmented Lagrangian can be
nonconvex, thus Lemma~\ref{thm:local_gradient_bound} cannot be invoked. However, the block coordinate descent algorithm is guaranteed to
converge (at least) to a stationary point as shown
in~\cite{xu2017globally}. This means that multiplier updates in ASYMM will surely occur after a finite number of primal minimizations. 
Moreover, as $\p^k$ grows, the Augmented Lagrangian typically becomes
locally strongly convex, see, e.g.,~\cite{bertsekas2014constrained}. If this
happens, the block coordinate descent algorithm approaches the corresponding
  minimum of the Augmented Lagrangian and, provided that for $k\geq\bar{k}$ the
local tolerances $\epsilon_i^k$ vanish as $k\to\infty$, the minimization of
the Augmented Lagrangian will be asymptotically exact.  This guarantees a
convergence result for ASYMM in the same sense as the one in the centralized
case: the convergence of the algorithm is contingent upon the
generation of (possibly local) minima of the Augmented Lagrangian that,
after some index $\hat{k}$, stay in the neighborhood of the same local minimum
$x^\star$ of problem~\eqref{pb:distributed_problem}.
As reported in~\cite{bertsekas2014constrained}, extensive numerical
  experience has shown that, from a practical point of view, choosing the obtained $\x^k$ as the initial condition for the $(k + 1)$-th
  minimization usually generates sequences $\{\x^k\}$ within a neighborhood of
  the same local minimum $\x^\star$.

\section{Numerical Results}
\label{sec:numerical}
Consider a network of $N$ sensors, deployed over a certain region, communicating
according to a connected graph $\GG=(\VV,\EE)$, which have to solve the optimization problem
\begin{equation*}
\begin{aligned}
& \m_x
& & \sum_{i=1}^N f_i(x)\\
& \st
& & \|x-c_i\|-R_i\leq 0, & i=1,...,N\\
& & & r_i-\|x-c_i\|\leq 0, & i=1,...,N,
\end{aligned}
\end{equation*}
which can be rewritten in the form of problem~\eqref{pb:distributed_problem}. 

Such a problem naturally arises, for example, in the context of source
localization under the assumption of unknown but bounded (UBB) noise, in which each agent knows its own absolute location $c_i$ and
takes a noisy measurement $y_i$ of its own distance from an an emitting source
located at an unknown location $x^\star$ as
$ y_i=\|x^\star-c_i\|+w_i$, with $|w_i|\leq \kappa_i$
for some $\kappa_i\geq 0$.

\begin{figure}
	\centering
	\begin{subfigure}{0.5\linewidth}
		\centering
		\includegraphics[width=\linewidth]{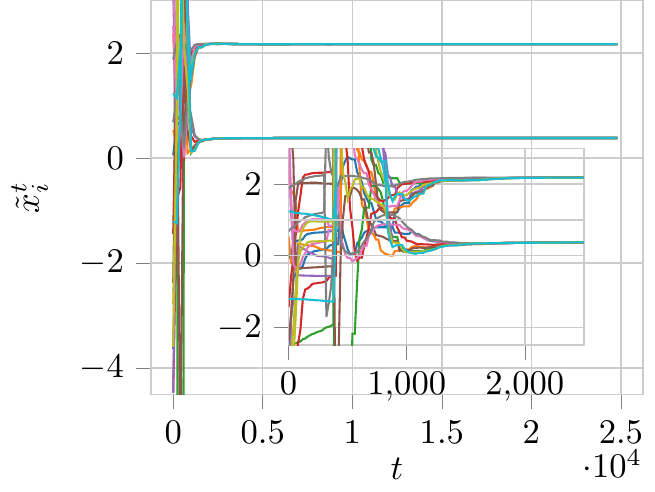}
		\caption{}
		\label{fig:xt}
	\end{subfigure}%
	\begin{subfigure}{0.5\linewidth}
		\centering
		\includegraphics[width=\linewidth]{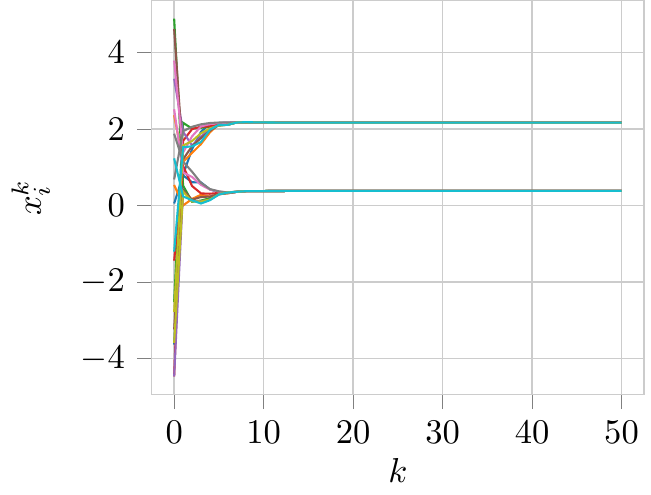}
		\caption{}
		\label{fig:xk}
	\end{subfigure}
	\caption{Evolution of the local decision variables $x_i$.}
	\label{fig:x}
\end{figure}
\begin{figure}
	\centering
	\includegraphics[width=0.55\linewidth]{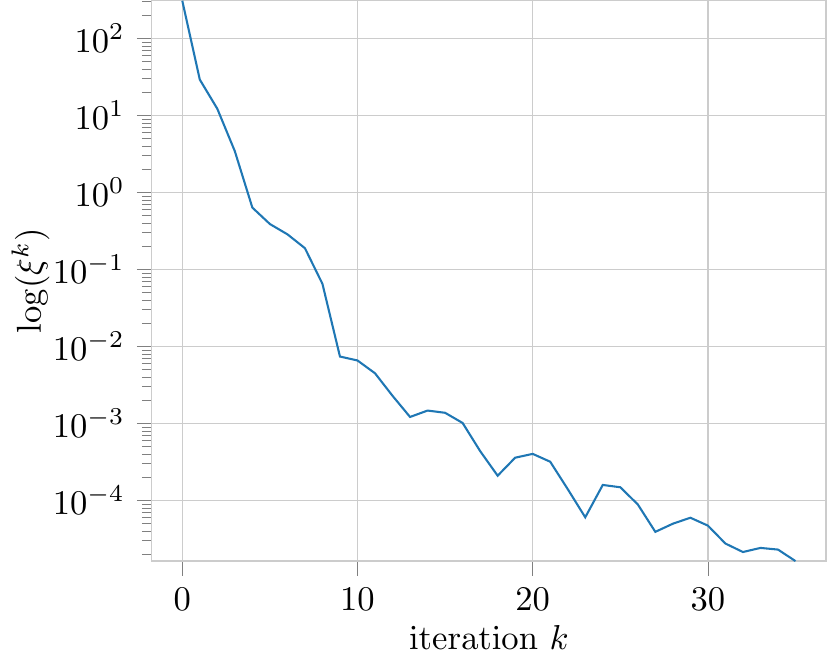}
	\caption{Logarithm of the measure of infeasibility $\xi^k$.}
	\label{fig:infeasibility}
\end{figure}

Suppose $f_i(x_i)=x_i^\top x_i$ for all $i\in\VV$. We report a simulation with
$N=10$ nodes and $n=2$, in which $x^\star\in U[-2.5,2.5]^n$,
$c_i\in U[-2.5,2.5]^n$ and $\kappa_i=U[0,0.3]$ for all $i\in\VV$, where $U[a,b]$
denotes the uniform distribution in the interval $[ a,b]$. The graph is modeled
through a connected Watts-Strogatz model in which nodes have mean degree $K=2$.
Let us define the measure of infeasibility at iteration $k$ as
$\xi^k=\sum_{i=1}^N[\max(0,\|x_i^k-c_i\|-R_i)+\max(0,r_i-\|x_i^k-c_i\|)+\sum_{j\in \NN_i\backslash i}\|x_i^k-x_j^k\|]$.
We run ASYMM for $25000$ iterations. In Fig.~\ref{fig:xt} the values $\tilde{x}_i^t$ are
shown for each $t=1,...,25000$ and each $i\in\VV$. A magnification of the first
$2500$ iterations is reported as a subplot. Fig.~\ref{fig:xk} shows the evolution
of $x_i^k$ for each $i\in\VV$. As it can be seen, the nodes performed $50$
multiplier updates along the $25000$ iterations. Finally, in Fig.~\ref{fig:infeasibility} the
values of $\xi^k$ are reported. 

\section{Conclusions}
\label{sec:conclusions}
In this paper we proposed ASYMM, an asynchronous distributed algorithm for a 
class of constrained optimization problems in which both local cost functions 
and constraints can be nonconvex. ASYMM has been proved to be equivalent to 
an inexact Method of Multipliers from which it inherits all the convergence 
properties. 


\begin{small}
	\bibliographystyle{IEEEtran}
	\bibliography{biblio}
\end{small}

\end{document}